\documentclass[a4paper,10pt,reqno]{amsart}
\usepackage{amscd}

\theoremstyle{plain} \numberwithin{equation}{section}
\newtheorem{thm}{Theorem}[section]

\newtheorem{prop}[thm]{Proposition}

\newtheorem{lemma}[thm]{Lemma}

\newtheorem{claim}[thm]{Claim}

\theoremstyle{definition}
\newtheorem{remark}[thm]{Remark}

\newtheorem{defn}[thm]{Definition}

\newtheorem{notn}[thm]{Notation}

\newtheorem{rmk}[thm]{Remark}

\newcommand{\bi}{\begin{itemize}}
\newcommand{\ei}{\end{itemize}}
\newcommand{\bp}{\begin{proof}}
\newcommand{\ep}{\end{proof}}

\def\AA{\mathbb{A}}
\def\CC{\mathbb{C}}
\def\GGG{\mathbb{G}}
\def\PP{\mathbb{P}}
\def\QQ{\mathbb{Q}}
\def\RR{\mathbb{R}}

\def\ZZ{\mathbb{Z}}

\def\ov{\overline}

\def\al{\alpha}
\def\be{\beta}

\def\ga{\gamma}

\def\la{\lambda}

\def\Ga{\Gamma}

\def\MM{\overline{M}}
\def\N{\mathcal{N}}
\def\O{\mathcal{O}}

\def\T{\mathcal{T}}

\def\W{\mathcal{W}}

\def\Cox{\text{\rm Cox}}
\def\PGL{\text{\rm PGL}}
\def\mult{\text{\rm mult}}

\def\hdeg{\text{hdeg}}

\def\Eff{\text{\rm Eff}}
\def\h{\text{h}}
\def\H{\text{H}}

\def\Hom{\text{Hom}}
\def\N{\text{N}}
\def\Pic{\text{\rm Pic}}
\def\Proj{\text{\rm Proj}}
\def\Spec{\text{Spec}}
\def\Bl{\text{\rm Bl}}

\def\dra{\dashrightarrow}
\def\hra{\hookrightarrow}

\def\ra{\rightarrow}

\def\bG{\Bbb G}
\def\bC{\Bbb C}
\def\g{{\mathfrak g}}
\def\ssl{{\mathfrak{sl}}}
\def\so{{\mathfrak{so}}}

\begin{document}

\title{Hilbert's 14-th Problem and Cox Rings}
\author{Ana-Maria Castravet and Jenia Tevelev}
\address{Department of Mathematics, University of Texas at Austin,
Austin, Texas, 78712}
\email{noni@math.utexas.edu and tevelev@math.utexas.edu}
\date{\today}

\begin{abstract}
Our main result is the description of generators of the total coordinate ring
of the blow-up of $\PP^n$ in any number of points that lie on a rational
normal curve.
As a corollary we show that the algebra of invariants of the action
of a two-dimensional vector group introduced by Nagata is finitely generated
by certain explicit determinants. We also prove the finite generation
of the algebras of invariants of actions of vector groups related
to T-shaped Dynkin diagrams introduced by Mukai.\end{abstract}

\maketitle


\section{Introduction}\label{intro}

Hilbert's 14-th Problem that we discuss is the following question: 
if an algebraic group $G$ acts linearly on a polynomial algebra $S$, is 
the algebra of invariants $S^G$ finitely generated? The answer is known to be
affirmative if $G$ is reductive (Hilbert \cite{Hi}) and if $G$
is the simplest nonreductive group $\bG_a$ (Weitzenb\"ock
\cite{Wei}). However, in general the answer is negative -- the
first counterexample was found by Nagata in 1958. Let
\begin{equation}\label{eqyationG}
G=\bG_a^g\subset \bG_a^r
\end{equation}
be a general linear subspace of codimension at least $3$. Consider
the following linear action of $\bG_a^r$ on
$S:=\bC[x_1,\ldots,x_r,y_1,\ldots,y_r]$: an element $(t_1,\ldots,
t_r)\in\bG_a^r$ acts by:
$$x_i\mapsto x_i,\quad y_i\mapsto y_i+t_i x_i, \quad 1\leq i\leq r.$$
The induced action of $G$ on $S$ is called {\em the Nagata
action}. The algebra  of invariants~$S^G$ is not finitely
generated if $g=13$ (Nagata \cite{Nag}), $g=6$ (Steinberg
\cite{St}), and finally $g=3$, $r=9$ (Mukai \cite{Mukai01}). Thus,
Hilbert's 14-th Problem has a negative answer for $\bG_a^3$. In
\cite{Mukai01}, Mukai asks what happens if $g=2$.

\begin{thm}\label{Thm_intro}
Assume without loss of generality
that $G=\bG_a^2\subset\bG_a^{n+3}$ is a linear subspace spanned by
rows of the matrix
$$\left[\begin{matrix}
1&1&\ldots&1\cr a_1&a_2&\ldots&a_{n+3}\cr
\end{matrix}\right]
$$
where $a_1,\ldots,a_{n+3}$ are general numbers. Then $S^G$ is
generated by $2^{n+2}$ invariants
\begin{equation}\label{invariants}
F_I=\left|\begin{matrix}
x_{i_1}&x_{i_2}&\ldots&x_{i_{2k+1}}\cr
a_{i_1}x_{i_1}&a_{i_2}x_{i_2}&\ldots&a_{i_{2k+1}}x_{i_{2k+1}}\cr
\vdots          &\vdots          &\ddots&\vdots\cr
a^{k}_{i_1}x_{i_1}&a^{k}_{i_2}x_{i_2}&\ldots&a^{k}_{i_{2k+1}}x_{i_{2k+1}}\cr
y_{i_1}&y_{i_2}&\ldots&y_{i_{2k+1}}\cr
a_{i_1}y_{i_1}&a_{i_2}y_{i_2}&\ldots&a_{i_{2k+1}}y_{i_{2k+1}}\cr
\vdots          &\vdots          &\ddots&\vdots\cr
a^{k-1}_{i_1}y_{i_1}&a^{k-1}_{i_2}y_{i_2}&\ldots&a^{k-1}_{i_{2k+1}}y_{i_{2k+1}}\cr
\end{matrix}\right|,
\end{equation}
where $I=\{i_1,\ldots,i_{2k+1}\}\subset\{1,\ldots,n+3\}$ is any
subset of odd cardinality $2k+1$.
\end{thm}

Of course,  it is possible that the algebra of invariants of
$\bG_a^2$ is not finitely generated for actions more complicated
than Nagata actions.

\

The ingenious insight of Nagata was to relate $S^G$ to a Cox ring.
Let $X$ be a projective algebraic variety over $\CC$. Assume that
divisors $D_1,\ldots,D_r$ freely generate the Picard group
$\Pic(X)$. Then the {\em Cox ring} of $X$ is the multigraded ring
$$\Cox(X)=\bigoplus_{(m_1,\ldots,m_r)\in\ZZ^{r}}
\H^0(X, m_1D_1+\ldots+m_r D_r)$$ (the basis is necessary to
introduce multiplication in a canonical way). This definition is a
generalization of the total coordinate ring of a toric variety
introduced by Cox \cite{Cox}. In fact, $\Cox(X)$ is isomorphic to
a polynomial ring if and only if $X$ is a toric variety
\cite[Prop.~2.10]{HuKeel}. For an arbitrary variety $X$, Hu and
Keel \cite[Prop.~2.9]{HuKeel} proved that $\Cox(X)$ is  finitely
generated if and only if $X$ is a \emph{Mori Dream Space}: 1) the
cone of nef divisors is generated by finitely many semi-ample line
bundles, and 2) the cone of moving divisors (divisors whose base
locus is of codimension at least $2$ in $X$) is the union of nef
cones of small modifications of $X$, i.e., varieties $X'$
isomorphic to $X$ in codimension $1$.

In the recent years, an explicit description of the ring $\Cox(X)$
has also proved useful for applications in arithmetic algebraic geometry.
Universal torsors were used for proving the Hasse Principle and weak
approximation for certain Del Pezzo surfaces
or for the counting of rational points of bounded height
(Colliot-Th\'el\`ene--Sansuc--Swinnerton-Dyer \cite{CS,CSS1,CSS2},
de la Bret\`eche \cite{Br}, Hassett--Tschinkel \cite{HT},
Salberger \cite{Sa}, Heath-Brown \cite{H-B}).

The relation to Nagata actions is as follows: if $G$ is as in
\eqref{eqyationG}, by \cite{Mukai01}, one has
$$S^G\simeq \Cox(\Bl_r\PP^{r-g-1})$$
where $\Bl_r\PP^{r-g-1}$ is the blow-up of $\PP^{r-g-1}$ at $r$
distinct points. Using this isomorphism, Theorem \ref{Thm_intro}
is equivalent to describing the Cox ring of a blow up of $\PP^n$
at $n+3$ points. It is a well-known fact that there is a unique
rational normal curve $C$ of degree $n$ in $\PP^n$ passing through
$n+3$ points in general position. We generalize Theorem
\ref{Thm_intro} as follows:

\begin{thm}\label{RationalNormalCurveTh}
Let $C\subset\PP^n$ be a rational normal curve of degree $n$ and
let $p_1,\ldots,p_r$ be distinct points on $C$, $r\ge n+3$. Let
$X=\Bl_{p_1,\ldots,p_r}\PP^n$. Then $\Cox(X)$ is finitely
generated by unique (up to scalar) global sections of exceptional
divisors $E_1,\ldots,E_r$ and divisors
\begin{equation}\label{MinimalDivisors}
E=kH-k\sum_{i\in I}E_i-(k-1)\sum_{i\in I^c}E_i
\end{equation}
for each subset $I\subset\{1,\ldots,r\}$, $|I|=n+2-2k$, $1\le k\le
1+n/2$. Here $H$ is the pull-back of the hyperplane class in
$\PP^n$.
\end{thm}

Geometrically, the divisors \eqref{MinimalDivisors} are proper
transforms of the following hypersurfaces in $\PP^n$ \cite[Ex.~9.6]{Harris}.
If $I$ is empty then \eqref{MinimalDivisors} is the $(n/2)$-secant
variety of~$C$. More
generally, if $\pi_I:\PP^n\dra\PP^{2k-2}$ is the projection from
the linear subspace spanned by the points $p_i$, $i\in I$ and
$C'=\pi_I(C)$, then $C'$ is a rational normal curve of degree
$2k-2$ and \eqref{MinimalDivisors} is the cone over
the $(k-1)$-secant variety of~$C'$.

An obvious generalization of Theorem~\ref{RationalNormalCurveTh}
would be to consider the Cox ring of the iterated blow up of $\PP^n$
along points, lines connecting them, $2$-planes, etc.
A special case of this construction is $\overline M_{0,n}$,
the Grothendieck--Knudsen moduli space of stable $n$-pointed rational curves.
If  the Cox ring of $\overline M_{0,n}$ is finitely generated, then
results of Keel--Hu \cite{HuKeel} and Keel--McKernan \cite{KM} almost
imply the ``Fulton conjecture'' for  $\overline M_{0,n}$
and therefore the description of the Mori cone of $\overline M_{g,n}$
(Gibney--Keel--Morrison \cite{GKM}).

\

Following Mukai \cite{Mukai02}, we also generalize Theorem
\ref{Thm_intro} in a different direction. Let $T_{a,b,c}$ be the
$T$-shaped tree with legs of length $a$, $b$, and $c$ with
$a+b+c-2$ vertices. We assume that $a,c\ge2$ and if $c=2$ then
$a>2$. Let
$$X_{a,b,c}=\Bl_{b+c}(\PP^{c-1})^{a-1}$$
be the blow-up of $(\PP^{c-1})^{a-1}$ in $r=b+c$ points in general
position. The effective cone $\Eff(X_{a,b,c})$ is the set of
effective divisors in $\Pic(X_{a,b,c})$. Mukai proves in
\cite{Mukai02} that if $T_{a,b,c}$ is not a Dynkin diagram of a
finite root system then $\Eff(X_{a,b,c})$ is not a finitely
generated semigroup and therefore $\Cox(X_{a,b,c})$ is not a
finitely generated algebra. Mukai also shows in \cite{Mukai02}
that the Cox algebra of any $X_{a,b,c}$ is isomorphic to the
algebra of invariants of a certain ``extended Nagata action''. We
deduce from Theorem~\ref{RationalNormalCurveTh}, using a trick from
commutative algebra, the following

\begin{thm}\label{MainTh}
The following statements are equivalent:
\begin{itemize}
\item $\Cox(X_{a,b,c})$ is a finitely generated algebra.
\item $\Eff(X_{a,b,c})$ is a finitely generated semigroup.
\item $T_{a,b,c}$ is a Dynkin diagram of a finite root system.
\item $\frac{1}{a}+\frac{1}{b}+\frac{1}{c}>1$.
\end{itemize}
Moreover, in these cases consider $Z=\Proj(\Cox(X))$ with respect
to the natural $\ZZ$-grading of $\Cox(X)$ defined in
\eqref{JHBJHGJGG}. Then $Z$ is a locally factorial,
Cohen--Macaulay, and Gorenstein scheme with rational
singularities. The Picard group $\Pic(Z)=\ZZ$ is generated by
$\O_Z(1)$ and the anti-canonical class is $-K_Z=\O_Z(d)$, where
$$d=abc\left(\frac{1}{a}+\frac{1}{b}+\frac{1}{c}-1\right)>0.$$
\end{thm}

The proof of the ``moreover'' part is exactly the same as Popov's
proof \cite{P} of the analogous statement for Del Pezzo surfaces
(or $X_{2,s-3,3}$ in our notation). We only sketch it for the
reader's convenience.

Explicitly, Theorem~\ref{MainTh} includes the following cases.
Mukai \cite{Mukai02} shows that $X_{a,b,c}$ is a small
modification of $X_{c,b,a}$, so we assume that $a\le c$ (if $X'$
is a small modification of $X$ then of course
$\Pic(X)\cong\Pic(X')$, $\Eff(X)\cong\Eff(X')$, and
$\Cox(X)\cong\Cox(X')$).

\begin{itemize}
\item $X_{2,2,n+1}=\Bl_{n+3}(\PP^n)$.

\item $X_{2,3,4}=\Bl_7\PP^3$, $X_{2,3,5}=\Bl_8\PP^4$.

\item $X_{3,2,3}=\Bl_5(\PP^2)^2$,  $X_{3,2,4}=\Bl_6(\PP^3)^2$,
$X_{3,2,5}=\Bl_7(\PP^4)^2$.

\item $X_{s+1,1,n+1}=\Bl_{n+2}(\PP^n)^s$. This case is well-known,
see Remark \ref{case r=n+2}.

\item Del Pezzo surfaces $X_{2,s-3,3}=\Bl_s\PP^2$, $s=4,5,6,7,8$.
In this case the finite generation of the Cox ring was proved by
Batyrev and Popov \cite{BP}.
\end{itemize}

We prove Theorems 1.1--1.3 in reverse order. In Section \ref{root_systems}
we describe the effective cone of $X_{a,b,c}$. In Section
\ref{finite_generation} we prove Theorem \ref{MainTh} (the finite
generation of $\Cox(X_{a,b,c})$) in all cases, except for $X_{2,3,4}$
and  $X_{2,3,5}$, for which the proof relies on the cases $n=3$ and $n=4$ of
Theorem \ref{RationalNormalCurveTh}. The latter is proved in Section
\ref{main}, which is the main section of the paper and is independent
of the previous sections. Theorem \ref{Thm_intro} is proved in the
concluding Section \ref{generators}. In particular, we prove the finite
generation of $\Cox(\Bl_{n+3}\PP^n)$ twice. First, we give a simple proof in
the framework of Theorem~\ref{MainTh}. Second, we give an independent proof
of the much stronger Theorem \ref{RationalNormalCurveTh} that gives explicit
generators for this ring. It is crucial for our proof to consider any number
of points on a rational normal curve. For example, finding generators for
$\Cox(\Bl_{n+3}\PP^n)$ relies on finding generators for the Cox ring of
the blow-up of $\PP^{n-1}$ in $n+3$ points lying on a rational normal curve,
etc, up to the blow-up of $\PP^2$ in $n+3$ points lying on a conic.
Our proof of Theorem~\ref{RationalNormalCurveTh} was inspired
by the ``whole-genome shotgun'' \cite{HG} method of genome sequencing
that involves breaking the genome up into very small pieces,
sequencing the pieces, and reassembling the pieces into the full genome
sequence. This method has some advantages (and disadvantages) over
the ``clone-by-clone'' approach that involves breaking the genome up into
relatively large chunks.

During the final stages of the preparation of this paper, Professor Shigeru
Mukai sent us his preprint \cite{Mukai05}, where he proves that the Cox ring of
$X_{2,b,c}$ is finitely generated when
$\frac{1}{2}+\frac{1}{b}+\frac{1}{c}>1$ by using a completely different
approach based on results of S. Bauer about parabolic bundles on curves.

{\bf Acknowledgements.} We would like to thank Professor Shigeru
Mukai for sending us his preprint \cite{Mukai05}. We are
especially grateful to Sean Keel for numerous discussions and
useful advice. We also thank Daniel Allcock, Gabi Farkas, Brendan
Hassett,  James McKernan, Dima Timashev, Nolan Wallach for useful
discussions. A-M. Castravet wishes to thank Institute des Hautes
{\'E}tudes Scientifiques for its hospitality during the summer of
2004.

\section{Root Systems and Effective Cones}\label{root_systems}

From now on we assume that $T_{a,b,c}$ is a Dynkin diagram of a
finite root system. It is well-known that this is equivalent to
$\frac{1}{a}+\frac{1}{b}+\frac{1}{c}>1$. Let
$$X=X_{a,b,c}.$$
The Picard group $\Pic(X)$ is a free $\ZZ$-module of rank
$a+b+c-1$ with a basis
$$H_1,\ldots,H_{a-1},\ \text{and}\ E_1,\ldots,E_r,$$
where $H_i$ is the pull-back of the hyperplane class from the
$i$-th factor of $(\PP^{c-1})^{a-1}$ and $E_j$ is the class of the
exceptional divisor over $p_j$, for $j=1,\ldots,r$, $r=b+c$.
We call this basis {\em tautological}.
If $a=2$ then we write $H$ instead of $H_1$ and make the appropriate
modifications in all notations. The anticanonical class of $X$ is
$$-K=c(H_1+\ldots+H_{a-1})-(ac-a-c)(E_1+\ldots+E_r).$$

Following \cite{Mukai02}, we define a symmetric bilinear form on
$\Pic(X)$ as follows:
\begin{equation}\label{bil_form}
(H_i,E_j)=0,\quad (H_i,H_j)=(c-1)-\delta_{i,j},
\quad (E_i,E_j)=-\delta_{i,j}.
\end{equation}

The following lemma is a straightforward calculation.

\begin{lemma}\cite{Mukai02}
$\Pic(X)$ has another $\ZZ$-basis $\al_1,\ldots,\al_{a+r-2},E_r$, where
$$\al_1=E_1-E_2,\ \ldots\ ,\ \al_{r-1}=E_{r-1}-E_r,$$
$$\al_r=H_1-E_1-\ldots-E_c,$$
$$\al_{r+1}=H_1-H_{2},\ \ldots\ ,\ \al_{a+r-2}=H_{a-2}-H_{a-1}.$$
Moreover, $\alpha_1,\ldots,\alpha_{a+r-2}$ is a $\ZZ$-basis of
the orthogonal complement $K^{\perp}$
and a system of simple roots of a finite root system with a
Dynkin diagram $T_{a,b,c}$.
\end{lemma}

Let $\W$ be the Weyl group generated by orthogonal reflections
with respect to $\alpha_1,\ldots,\al_{a+r-2}$. Then $K$ is
$\W$-invariant. Mukai calls $D\subset X$ a $(-1)$-divisor if there
is a small modification $X\dra X'$ such that $D$ is the
exceptional divisor for a blow-up $X'\ra Y$ at a smooth point.
Note that any $(-1)$-divisor must appear in any set of generators
of $\Eff(X)$.

\begin{lemma}\cite{Mukai02}\label{cremona}
For each transformation $w:\,\Pic(X)\to\Pic(X)$ of $\W$, there is a small
modification $X\dra X_w$ with the following property. $X_w$ is
also a blow-up of $(\PP^{c-1})^{a-1}$ in $r=b+c$ points
$q_1,\ldots,q_r$ in general position and the pull-back of the
tautological basis of~$X_w$ coincides with the transformation of the
tautological basis of $X$ by $w$. In particular, every divisor
$E\in\W\cdot E_r$ is a $(-1)$-divisor and $H^0(X,E)$ is spanned by
a single section $x_E$.
\end{lemma}

The proof is an application of Cremona transformations.
The case $a=2$ appeared in~\cite{Do2}
(where it is attributed to Coble).
The case $a=2$, $c=3$ is well-known
from the theory of marked Del Pezzo surfaces.

\begin{lemma}\label{preserve}
The action of $\W$ on $\Pic(X)$ preserves $\Eff(X)$.
\end{lemma}

\bp Let $D\in\Pic(X)$ and $w\in\W$. We claim that $H^0(X,D)\simeq
H^0(X,w\cdot D)$. We have
$$D=d_1H_1+\ldots+d_{a-1}H_{a-1}-m_1E_1-\ldots-m_rE_r.$$
$H^0(X,D)$ can be identified with the subspace of polynomial
functions on $(\CC^c)^{a-1}$ of multidegree $(d_1,\ldots,d_{a-1})$
vanishing to the order at most $m_i$ at the point $p_i$. By
Lemma~\ref{cremona}, $H^0(X,w\cdot D)$ has the same interpretation
for another choice of general points  $q_1,\ldots,q_r$. Now the
claim follows from semi-continuity if the points $p_1,\ldots,p_r$
are sufficiently general. \ep

Let $\Eff_\RR(X)\subset\Pic(X)\otimes\RR$ be the cone spanned by
$\Eff(X)$. Let $\N_1(X)$ be the group generated over $\ZZ$ by
$1$-cycles on $X$ modulo rational equivalence. Intersection of
cycles gives a non-degenerate pairing $\Pic(X)\times
\N_1(X)\ra\ZZ$.  For $i=1,\ldots,a-1$, let $l_i\in\N_1(X)$ be the
class of the proper transform of a general line in the $i$-th copy
of $\PP^{c-1}$. For $i=1,\ldots,r$, let $e_i\in\N_1(X)$ be the
class of a general line in $E_i$. Then it is easy to check that

\begin{equation}\label{pairing}
H_i\cdot l_j=\delta_{i,j},\quad
H_i\cdot e_j=0,\quad
E_i\cdot e_j=\delta_{i,j}.
\end{equation}

Since the intersection pairing is non-degenerate,
it follows that $\N_1(X)$ is generated over $\ZZ$ by the
classes $l_1,\ldots,l_{a-1}, e_1,\ldots, e_r$.
The action of $\W$ on $\Pic(X)$ induces an action on $\N_1(X)$.

A class $\ga$ in $\N_1(X)$ is called \emph{nef} if for any
effective divisor $D$ on $X$, $D\cdot\ga\geq0$.

\begin{lemma} \label{nef_classes}
The classes $l_i$, $l_1+\ldots+l_{a-1}-e_i$ are nef, for all $i=1,\ldots,r$.
\end{lemma}
\bp  Note that if a family of curves with class $f$ covers $X$
(i.e., through a general point of $X$ there is an irreducible curve in the
family that passes through it), then $f$ is a nef class: if $D$ is
an effective divisor, there is an irreducible curve in the family that is not
contained in $D$, therefore, $D\cdot f\ge0$. This is obviously the
case if $f=l_i$. If $f=l_1+\ldots+l_{a-1}-e_i$, then $f$ is the
proper transform in $X$ of a curve of multidegree $(1,\ldots,1)$
in $(\PP^{c-1})^{a-1}$ that passes through the point $p_i$.
This family contains an irreducible curve by
Bertini's theorem and we can use the $2$-transitive action of
$(\PGL_c)^{a-1}$ on $(\PP^{c-1})^{a-1}$ to find a curve through any point.\ep

\begin{defn} Define the {\em degree} of $D\in\Pic(X)$
as an integer $$\deg(D)=\frac{1}{ac-a-c}(D, -K).$$
Clearly, $\deg D$ is $\W$-invariant
and any divisor in the orbit $\W.E_r$ has degree $1$.
\end{defn}

\begin{defn}
Let $\g_{a,b,c}$ be a semisimple Lie algebra with the Dynkin
diagram $T_{a,b,c}$. Let $\Lambda\subset K^\perp\otimes\QQ$ be the
{\em weight lattice} spanned by {\em fundamental weights}
$\omega_1,\ldots,\omega_{a+r-2}$ defined by
$(\omega_i,\alpha_j)=\delta_{i,j}$. For any $\omega\in\Lambda$,
let $L_\omega$ be an irreducible $\g_{a,b,c}$-module with the
highest weight $\omega$, see for example \cite{VO}. $L_\omega$ is
called {\em minuscule} if weights $\W\cdot\omega$ are its only
weights. Let $\pi:\,\Pic(X)\to K^\perp\otimes\QQ$ denote the
orthogonal projection.
\end{defn}

\begin{thm}\label{effective_cone}\label{deg1_divisors}
$\Eff(X)$ is generated as a semigroup by divisors of degree~$1$.
$\Eff_\RR(X)$ is generated as a cone by $D\in\W\cdot E_r$.
Projection $\pi$ induces a bijection between divisors of degree~$1$
and weights of $L_{\omega_{r-1}}$ such that divisors
in $\W\cdot E_r$ correspond to weights in $\W\cdot\omega_{r-1}$.
In particular, $L_{\omega_{r-1}}$ is minuscule if and only if
the only effective divisors of degree $1$ are $D\in\W\cdot E_r$.
\end{thm}

\begin{rmk}\label{kljhlkjglkg}
The classification of minuscule representations is well-known.
The only arising cases are
$$\Bl_{n+3}\PP^n,\quad \Bl_{n+2}(\PP^n)^s,\quad
\hbox{\rm and}\quad \Bl_s\PP^2\ (s=4,5,6,7).$$
If  $X=\Bl_{n+3}\PP^n$ then  $L_{\omega_{r-1}}$ is a halfspinor
representation of $\so_{2n+6}$.
Here's another example: let $X=X_{2,3,3}$ be the blow-up of $\PP^2$ in
$6$ general points,
i.e.~a smooth cubic surface.
Divisors of degree $1$ are the $27$ lines.
The corresponding minuscule representation
$L_{\omega_{r-1}}$ is the $27$-dimensional representation
of $E_6$ as a Lie algebra of infinitesimal norm similarities of the
exceptional Jordan algebra.
\end{rmk}

\bp[Proof of Theorem~\ref{effective_cone}]
Let $\Ga_k$ be the intersection of the convex hull of
$W\cdot(kE_r)$ with $\Eff(X)$ and let
$\Ga\subset\Pic(X)\otimes\RR$ be the cone spanned by $\Ga_1$.
Since $\pi(E_r)=\omega_{r-1}$ and any element of $K^{\perp}$ is an
integral combination of roots, it follows from the basic
representation theory of semisimple Lie algebras \cite{VO} that
$\pi(\Ga_k)$ is the set of weights of an irreducible $\g$-module
$L_{k\omega_{r-1}}$ with the highest weight $k\omega_{r-1}$.
Since $L_{k\omega_{r-1}}\subset L_{\omega_{r-1}}^{\otimes k}$
($L_{k\omega_{r-1}}$ is the so-called {\em Cartan component} of
$L_{\omega_{r-1}}^{\otimes k}$), any weight in
$\pi(\Ga_k)$ is a sum of $k$ weights from $\pi(\Ga_1)$, and
therefore any divisor in $\Ga_k$ is a sum of $k$ divisors
from~$\Ga_1$. It follows that $\Eff(X)\cap\Ga$ is generated by
$\Ga_1$ as a semigroup.

It remains to show that $\Eff(X)_{\RR}\subset\Ga$.
We will find all faces of $\Ga$ and show that the inequalities
that define them are satisfied by any effective divisor.

By Lemma \ref{preserve}, it suffices to find faces of $\Ga$
adjacent to the ray spanned by $E_r$ up to the action of the
stabilizer of $E_r$ in $\W$. The algorithm for finding faces of
these so-called Coxeter polytopes is explained, for example, in
\cite{C}, p. 9. They are in one-to-one correspondence with
connected maximal subdiagrams of $T_{a,b,c}$ that contain the
support of the highest weight, i.e., ~the node that corresponds to
the simple root $\alpha_{r-1}$ in our case. There are two types of
such diagrams given by roots:
$$1)\ \al_2,\al_3,\ldots,
\al_{a+r-2}\quad\hbox{\rm and}\quad 2)\ \al_1, \al_2,\ldots, \al_{a+r-3}.$$
For each subdiagram,
the linear span of the corresponding face
is spanned by simple roots in the subdiagram and by $E_r$.

Using formulas \eqref{pairing},
any face of $\Ga$ is given (up to the action of $\W$) by inequality
\begin{equation}\label{inequality}
D\cdot f\ge0,
\end{equation}
where
$$1)\ f=l_1+\ldots+l_{a-1}-e_1
\quad\hbox{\rm and}\quad 2)\ f=l_{a-1}.$$
By Lemma \ref{nef_classes}, the
class $f$ is nef. Hence, for any $D$ effective, $D\cdot f\ge0$. We
conclude that \eqref{inequality} is, in fact, satisfied by  any
effective divisor and $\Eff_\RR(X)=\Ga$. \ep

\section{Proof of Theorem \ref{MainTh}}\label{finite_generation}

The following is a direct generalization from \cite{BP} (Prop. 4.4):

\begin{prop}\label{charts}
Let $\pi:\,X\to X'$ be the blow up of a smooth point.
Let $E\subset X$ be an exceptional
divisor, and let $x_E\in H^0(X,E)\subset\Cox(X)$ be the
corresponding section. Then there is an isomorphism of rings
$$\Cox(X)_{x_E}\cong \Cox(X')[T,T^{-1}].$$
\end{prop}

\bp
Any divisor $D\in\Pic(X)$ can be uniquely written as
$D=D_0-mE$, where $D_0\in\pi^*\Pic(X')$, $m\in\ZZ$.
We identify $\Pic(X')$ with $\pi^*\Pic(X')\subset\Pic(X)$
and $\Cox(X')$ with $\pi^*\Cox(X')\subset\Cox(X)$.
The latter embedding extends to a ring homomorphism
$$\Cox(X')[T,T^{-1}]\ra \Cox(X)_{x_E}$$
by sending $T$ to $x_E$. We show that this is an isomorphism by
constructing an inverse to it.

If $m\ge0$ and $s$ is a section in $\H^0(X, D)$, then let
$s_0=s\cdot x_E^m\in\H^0(X', D_0)$. Define a map
$$\H^0(X,D)\ra\H^0(X',D_0)T^{-m},\quad s\mapsto s_0T^{-m}.$$

If $m<0$ then the canonical inclusion $\H^0(X,D_0)\hra\H^0(X,D)$
is an isomorphism. To see this, note that for any $i\geq0$ there
is an exact sequence
$$0\ra\H^0(X, D_0+iE)\ra\H^0(X,D_0+(i+1)E)\ra
\H^0(E,(D_0+(i+1)E)|_E)=0,$$
where the last equality follows from
$$\O(D_0)|_{E}=\O_E\quad\hbox{\rm and}\quad \O(E)|_{E}=\O_E(-1).$$
Define a map $\H^0(X,D)\ra\H^0(X',D_0)T^{-m}$ in the same way, by
sending $s$ to $s_0\cdot T^{-m}$, where $s\in\H^0(X, D)$ is the
image of a section $s_0\in\H^0(X,D_0)$. This gives a map
$\Cox(X)\ra \Cox(X')[T, T^{-1}]$ which maps $x_E$ to $T$. One can
check directly that this is a ring homomorphism. The induced map
$\Cox(X)_{x_E}\ra \Cox(X')[T, T^{-1}]$ is the desired inverse. \ep

\begin{notn}
In this section,
$$X=X_{a,b,c}.$$
\end{notn}
\begin{prop}\label{CoxUFD}
$\Cox(X)$ is a unique factorization domain.
\end{prop}

\bp The Cox ring of a normal projective variety is known to be a
UFD \cite{EKW}. We can also use a simple observation: the ring of
invariants of a UFD with respect to the action of a connected
algebraic group without nontrivial characters is a UFD, see
\cite{PV}. By \cite{Mukai02}, $\Cox(X)$ is a ring of invariants of
an {\em extended Nagata action}. \ep

\begin{defn}\label{JHBJHGJGG}
We define a $\ZZ$-grading of $\Cox(X)$ by $\deg(s)=\deg(D)$ for any
$s\in\H^0(X, D)$. In particular, $\deg(x_E)=1$ for any $E\in \W\cdot E_r$.
\end{defn}

\begin{defn}
Let $\Cox'(X)\subset\Cox(X)$ be a subalgebra generated by sections
$x_E$, for $E\in\W\cdot E_r$. We say that $\Cox(X)$ is
\emph{minuscule} if $\Cox(X)=\Cox'(X)$.
\end{defn}

\begin{defn}Let $\PP(X)=\Proj(\Cox(X))$, $\AA(X)=\Spec(\Cox(X))$, and
$Z=\Proj(\Cox'(X))$, where $\Cox(X)$ and $\Cox'(X)$ are considered
with their $\ZZ$-grading as in Definition~\ref{JHBJHGJGG}
(we will show that in fact $Z\cong\PP(X)$).
\end{defn}

Inspecting the list of all possible $X_{a,b,c}$ given in the
introduction, we see that $X_{a,b-1,c}$ is contained in the
following list:

\begin{itemize}
\item $X_{s+1,1,n+1}=\Bl_{n+2}(\PP^n)^s$. This variety is minuscule,
see Remark~\ref{case r=n+2}.
\item Del Pezzo surfaces $X_{2,s-3,3}=\Bl_s\PP^2$,
$s=4,5,6,7$. In this case $\Cox(X)$ is minuscule by a theorem of
Batyrev and Popov \cite{BP}.
\item $X_{2,2,4}=\Bl_6(\PP^3)$,
$X_{2,2,5}=\Bl_7(\PP^4)$. These varieties are also minuscule
by our Theorem~\ref{Thm_intro} (which will be proved later).
\end{itemize}
Therefore, $X_{a,b-1,c}$ is minuscule in all cases.

Let $R=\Cox(X)$, $R'=\Cox'(X)$, $R_0=\Cox(X_{a,b-1,c})$. Let $Q$
be the field of fractions of $R$. We claim that $R$ is contained
in all the localizations $R'_{x_E}\subset Q$. By Lemma
\ref{cremona}, there is a small modification $\tilde X$ of $X$
isomorphic to $\Bl_r(\PP^{c-1})^{a-1}$, the blow-up of
$(\PP^{c-1})^{a-1}$ in $r=b+c$ points $q_1,\ldots,q_r$ in general
position, such that the pullback of $E$ is contracted to $q_r$. By
Proposition \ref{charts}, $R\subset (R_0)_{x_E}$. It remains to
notice that $R_0\subset R'$ because $R_0$ is minuscule.

\begin{claim}
$R$ is integral over $R'$
\end{claim}

\bp This is a standard proof, see for example \cite[p. 123]{Ha}.
Let $z\in R$ be a homogeneous element of a positive degree.
To show that $z$ is integral over $R'$, it suffices to find a faithful
$R'[z]$-module $M$ finitely generated as an $R'$-module.
Let $M$ be the set of elements in $R'$ of degree greater than $N$,
where $N$ has to be chosen adequately.
Obviously, $M$ is an $R'[z]$-module if $zM\subset R'$.
So choose $N$ to be $kn+1$, where $k$ is the number of generators $x_i$
in $R'$, and $n$ is the maximum of integers $n_i$ such that
$zx_i^{n_i}\in R'$. Clearly, $M$ is a finitely generated $R'$-module.
Since $R$ is a domain, $M$ is of course a faithful $R'[z]$-module.
\ep

It follows that $R$ is
integral over $R'$ and, therefore, $R$ is finitely generated.

Now we prove the ``moreover'' part of the theorem following
Popov's proof \cite{P} of the analogous statement for Del Pezzo
surfaces.

For each $E\in\W\cdot E_r$ consider the open chart $U_E(X)\subset
Z$ given by $x_E\neq0$. These charts cover $Z$.
Let $U_E'(X)\subset\PP(X)$ be a chart given by $x_E\neq0$.
Since $R$ is integral over $R'$, it is easy to see
that the radical of the ideal of $R$ generated by $x_E$'s
is the irrelevant ideal. It follows that charts $U_E'(X)$
cover $\PP(X)$. Since $R\subset\cap_ER'_{x_E}$, we have
$R'_{x_E}= R_{x_E}$ for any
$E\in\W\cdot E_r$.
It follows that, in fact, $U_E(X)\simeq U_E'(X)$,
the inclusion $R'\subset R$ induces an isomorphism $\phi:\PP(X)\ra Z$,
and $\phi^*\O_Z(m)\cong \O_{\PP(X)}(m)$. Moreover, it is true in general 
that if a graded ring $R$ is a UFD and the irrelevant ideal is the radical 
of the ideal generated by degree $1$ elements, then the Picard group 
of $\Proj(R)$ is $\ZZ$ and it is generated by $\O(1)$.

It follows from
Proposition~\ref{charts} that $U_E(X)\cong \AA(X_{a,b-1,c})$ is
factorial by Proposition~\ref{CoxUFD}. Therefore, $Z$ is locally
factorial and, in particular, $Z$ is normal.

Arguing by induction on $b$,
we can assume that all statements of Theorem \ref{MainTh} are satisfied for
$Y=X_{a,b-1,c}$.
Let $W=\PP(Y)$. Thus $W$ is a Cohen--Macaulay and Gorenstein scheme
with rational singularities, $\Pic(W)=\ZZ$ is generated by
$\O_W(1)$ and the anti-canonical line bundle $\omega_W$ is ample.

\begin{lemma}
$H^i(W,\O(k))=0$ for $i\ge1$, $k\ge0$.
\end{lemma}

\bp Notice that $\O(k)=\omega_W\otimes L$ with $L$ ample. Let
$\pi:\,\tilde W\to W$ be a resolution of singularities. Then
$H^i(\tilde W,\omega_{\tilde W}\otimes\pi^*(L))=0$ by Kodaira
vanishing because $\pi^*(L)$ is big and nef. Now use the Leray
spectral sequence and the definition of rational singularities
($R^i\pi_*\omega_{\tilde W}=0$ for $i>0$) to conclude that
$H^i(W,\O(k))=0$. \ep

Since $Y$ is minuscule, $W$ is projectively normal in the
projective embedding given by $\O_{W}(1)$. Note that
$U_E(X)\cong\AA(Y)$ is an affine cone over $Y$.
It follows that $\AA(Y)$ has rational singularities
by~\cite[Theorem 1]{KR} and therefore is Cohen--Macaulay
\cite{Ke}. Since $\AA(Y)$ is factorial and Cohen--Macaulay, it is
Gorenstein \cite[Ex. 21.21]{Ei}.

It remains to calculate the anticanonical class of $\PP(X)$. By
\cite{HuKeel}, $X$ is the GIT quotient of $\AA(X)$ for the action
of the torus $\Hom(\Pic(X),\GGG_m)=\GGG_m^{r+1}$. Moreover, $X$ is
the GIT quotient of $\PP(X)$ for the induced action of $\GGG_m^r$.
Let $U$ be the semi-stable locus in $\PP(X)$. Note that there are
no strictly semi-stable points \cite[Prop. 2.9]{HuKeel}.  It is
easy to see by induction using charts $U_E(X)$ that $\GGG_m^r$
acts on $\PP(X)$ with connected stabilizers. By Luna's \'etale
slice theorem \cite[p. 199]{Mumford}, this implies that
$\pi:\,U\to X$ is a principal \'etale fibre bundle. In particular,
$U$ is smooth. By the general theory of Cox varieties \cite[Prop. 
2.9]{HuKeel}, $\PP(X)\setminus U$ has codimension at least $2$ in
$\PP(X)$, and therefore $\Pic(U)\cong\ZZ\{\O(1)\}$. By the GIT,
the pull-back map $\pi^*$ between the Picard groups is the map
given by degree: $\pi^*(D)=\deg(D)$.

It is enough to prove that $K_U=\O_U(-d)$. Let $\T_X$ (resp.
$\T_U$) be the tangent sheaf of $X$ (resp. $U$). There is an exact
sequence of locally free sheaves:
$$0\ra \O_U^r\ra\T_U\ra\pi^*\T_X\ra0$$
(the relative tangent sheaf of a principal \'etale bundle is
canonically a trivial bundle with fiber isomorphic to the Lie
algebra of $\GGG_m^r$). Taking Chern classes, it follows that
$c_1(\T_U)=\pi^*(c_1(\T_X))$; hence, $-K_U=\pi^*(-K_X))=\O(d)$,
where
$$d=\deg(-K_X)=abc\left(\frac{1}{a}+\frac{1}{b}+\frac{1}{c}-1\right).$$
\qed

\begin{remark}\label{case r=n+2}
Here we consider the case of $X=X_{s+1,1,n+1}=\Bl_{n+2}(\PP^n)^s$.
Then it is well-known
and easy to check that $X$ is the GIT quotient of the Grassmannian
$G(s+1, n+s+2)$. It follows from \cite{HuKeel} that $\Cox(X)$ is
isomorphic to the total coordinate ring of $G(s+1, n+s+2)$ which is
generated by the $n+s+2\choose s+1$ Pl\"ucker coordinates. On the
other hand, the orbit $\W\cdot E_r$ in this case consists of
precisely $n+s+2\choose s+1$ divisors, the dimension of the
minuscule representation of $\g_{s+1,1,n+1}=\ssl_{n+s+2}$ in
$L_{\omega_n}=\Lambda^{s+1}\CC^{n+s+2}$. It follows that $\Cox(X)$
is minuscule.
\end{remark}

\section{Proof of Theorem~\ref{RationalNormalCurveTh}}\label{main}

\begin{notn}
Let $X=\Bl_r\PP^n$ be the blow-up of $\PP^n$ at $r$ distinct points
$p_1,\ldots,p_r$ ($r\geq n+3$) that lie on a rational normal curve $C$
of degree $n$. Let $E_1,\ldots,E_r$ be the
exceptional divisors and $H$ the hyperplane class.
Let
$$\al=r-n-2.$$
Let $\tilde C$ be the proper transform of $C$ on $X$.
\end{notn}

\begin{lemma} \label{lower_bound_multip}
Let $D\subset\PP^n$ be a hypersurface of degree $d$ that contains
$C$ with multiplicity $m$. If $D$ has multiplicity $m_i$ at $p_i$,
$i=1,\ldots,r$, then one has:
$$m\geq\frac{\sum_{i=1}^r m_i-nd}{\al}.$$
\end{lemma}

\bp
Recall that the multiplicity of a divisor along a curve is the
multiplicity at a general point of a curve.
Let $\tilde{D}$ be the proper
transform of $D$ on $X$.
Let $\pi':X'\ra X$ be
the blow-up of $X$ along~$\tilde{C}$ and let $E$ be the
exceptional divisor. Then $E\cong\PP(N_{\tilde{C}|X})$, where
$N_{\tilde{C}|X}$ is the normal bundle of $\tilde{C}$ in $X$.
One has
$$N_{C|\PP^n}\cong \O(n+2)^{\oplus (n-1)},$$ and therefore
$$N_{\tilde{C}|X}\cong\pi^*N_{C|\PP^n}\otimes \O_X(-E_1-\ldots-E_r) \cong
\O(n+2)^{\oplus (n-1)}\otimes \O(-r)\cong
\O(-\al)^{\oplus (n-1)}.$$

It follows that $E\cong\PP^1\times\PP^{n-2}$.
Let
$$q_1:\PP^1\times\PP^{n-2}\ra\PP^1,\quad
q_2:\PP^1\times\PP^{n-2}\ra\PP^{n-2}$$
be the two projections.
Then $\O(E)_{|E}\cong q_1^*\O(-\al)\otimes q_2^*\O(-1)$.

Let $D'$ be the proper transform of $\tilde{D}$ on $X'$. Then
${\pi'}^*\tilde{D}=D'+m E$. Denote
$$a=-\tilde{D}.\tilde{C}=\sum_{i=1}^rm_i-nd.$$

Note that ${\pi'}^*{\O_X(\tilde D)}_{|E}=q_1^*\O(-a)$. Since
${\O_{X'}(D')}_{|E}=q_1^*\O(-a+m\al)\otimes q_2^*\O(m)$ is an
effective divisor on $E$, it follows that $-a+m\al\geq0$. Hence,
$m\geq a/\al$. \ep

\begin{lemma}\label{minimal}
Consider the divisor~\eqref{MinimalDivisors} on $X$.
Then $E$ is the
proper transform of a unique hypersurface of degree $k$ in $\PP^n$ that
has multiplicity $k$ at any $p_i$ with $i\in I$ and $k-1$ at
all other points of $C$. In particular, $\H^0(X,E)\cong\CC$ and
$E-E_i$ is not effective for any $i=1,\ldots,r$.
\end{lemma}

\bp Let $J\subset I^c$ be any subset with $|J|=2k+1$. The divisor
$$E'=kH-k\sum_{i\in I}E_i-(k-1)\sum_{i\in J}E_i$$
is an effective divisor of degree $1$ on the blow-up
$\Bl_{n+3}\PP^n$ of $\PP^n$ along the points~$p_i$ for $i\in I\cup
J$. It follows that $\h^0(X,E')=1$ and, for any $i\in I\cup J$,
the divisor $E-E_i$ is not effective. It follows that $E$ is the
proper transform of a unique hypersurface $Z$ of degree~$k$ in $\PP^n$
such that
$$\mult_{p_i}Z=k\quad (i\in I)\quad\hbox{\rm and}\quad
\mult_{p_i}Z=k-1\quad (i\in J).$$
Since $Z$ is the image of $E$, and therefore does not depend on the choice of $J$,
we have $\mult_{p_i}Z=k-1$ for any $i\in I^c$.
If $p$ is a
point on $C$ different than $p_1,\ldots,p_r$, consider the variety
$\Bl_{r+1}\PP^n$ that is the blow-up of $X$ at $p$. Let $E_{r+1}$
be the exceptional divisor. By applying the same argument to the
divisor $E-(k-1)E_{r+1}$ on $\Bl_{r+1}\PP^n$, it follows that the
multiplicity of $Z$ at $p$ is exactly $k-1$. \ep

\begin{defn} We call the divisors $E$ in \eqref{MinimalDivisors}
\emph{minimal} divisors on $\Bl_r\PP^n$.
We call an element in $\Cox(X)$ a \emph{distinguished section} if
it is a monomial in the sections $x_E\in\H^0(X,E)$, where $E$ is
either a minimal divisor on $X$ or an exceptional divisor~$E_i$.
The ring $\Cox(X)$ is minuscule if it is generated by
distinguished sections.
\end{defn}

We prove that
$\Cox(X)$ is minuscule by induction on
$n$ and $r$. Theorem \ref{n=2}  proves this for $n=2$.
Assume from now on that $n\geq3$.

\begin{defn}
Let
\begin{equation}\label{divisorD}
D=dH-\sum_{i=1}^rm_i E_i
\end{equation}
be any divisor on $X$.
We call $d$ the {\em $H$-degree} of $D$, denoted by $\hdeg(D)$.
\end{defn}

\begin{notn}
Consider the projection $\pi_1:\PP^n\dra\PP^{n-1}$ from $p_1$
and let $q_i=\pi(p_i)$ for $i=2,\ldots,r$.
Note that $q_2,\ldots,q_r$ lie on a rational normal curve $\pi_1(C)$ of
degree $n-1$ in~$\PP^{n-1}$.
Let
$Y=\Bl_{r-1}\PP^{n-1}$ be the blow-up of $\PP^{n-1}$ at
$q_2,\ldots,q_r$.
Let $\ov E_2,\ldots,\ov E_r$ be the
exceptional divisors on $Y$ and $\ov H$ the hyperplane class.
Consider the linear map $\Pic(X)\to\Pic(Y)$ that maps \eqref{divisorD} to
\begin{equation}\label{res_divisor}
\tilde{D}= m_1\ov H-\sum_{i=2}^r(m_i+m_1-d)\ov E_i.
\end{equation}
\end{notn}

\begin{lemma}\label{trick}
If $\hdeg(D)=\hdeg(D')$ and $\tilde{D}=\tilde{D'}$ then $D=D'$.
\end{lemma}

\bp
This is because by (\ref{res_divisor}), $\tilde{\Delta}=0$ implies that
$\Delta=e(H-\sum_{i=1}^rE_i)$, for some $e\in\ZZ$. Hence, if the
$H$-degree of $\Delta$ is $0$, then $\Delta=0$.
\ep

\begin{lemma}
There is a map $r$
that makes the following diagram commutative:
$$
\begin{CD}\label{res2}
\H^0(X,D)@>r>>\H^0(Y,\tilde{D})\\
@Vr'VV            @VViV \\
\H^0(E_1,D_{|E_1})@=\H^0(\PP^{n-1},\O(m_1))\\
\end{CD}
$$
Here $r'$ is the restriction map and $i$
is the canonical injective map
given by push-forward.
For any divisors $D_1$,
$D_2$ on $X$ and $s_1\in\H^0(X,D_1)$, $s_2\in\H^0(X,D_2)$, if
$D=D_1+D_2$, then
$$\tilde{D}=\tilde{D_1}+\tilde{D_2},\quad
r(s_1s_2)=r(s_1)r(s_2).$$
\end{lemma}

\bp We can identify $E_1$ with the image of the projection $\pi_1$
and view $r'$ as a map
$$r':\,\H^0(X,D)\to\H^0(\PP^{n-1},\O(m_1))=\H^0(Y,m_1\ov H).$$
Let $x_{\ov E_i}$ be a generator for $\H^0(Y,\ov E_i)\cong\CC$.
Note
that if for some $i=2,\ldots,r$ one has $m_1+m_i-d>0$,
then the image of $r'$ lies in the linear
subsystem
$$|m_1\ov H-(m_1+m_i-d)\ov E_i|\subset|m_1\ov H|$$
and therefore
$r'(s)$ is
divisible by $x_{\ov E_i}^{-d+m_1+m_i}$
for any $s\in\H^0(X,D)$.
It follows that we can formally define
$$r(s)=r'(s)\prod_{i=2}^rx_{\ov E_i}^{d-m_1-m_i}.$$
The last statement of the lemma is clear.
\ep

\begin{rmk} The geometric interpretation for the map $r$ is
as follows.
Let $l_{i,j}$ be the
proper transform on $X$ of the line in $\PP^n$ joining the points
$p_i$ and $p_j$. Then $q_2,\ldots,q_r$ are the points on
$E_1\cong\PP^{n-1}$ where $l_{1,2}$, $l_{1,3}$, $\ldots$, $l_{1,n}$
intersect $E_1$.
Let $\tilde{X}$ be the blow up of $X$ along
$l_{1,2},\ldots,l_{1,n}$ and let $E_{1,2},\ldots, E_{1,n}$ be the
exceptional divisors. The normal bundle $N_{l_{i,j}|X}$ of
$l_{i,j}\cong\PP^1$ in $X$ is $\O(-1)^{\oplus (n-1)}$. The
exceptional divisors $E_{1,j}$ are given by:
$$E_{1,i}\cong\PP(N_{l_{i,j}|X})\cong
l_{i,j}\times\PP^{n-2}\cong\PP^1\times\PP^{n-2}.$$

For any $n\geq3$, there is morphism $\tilde{X}\ra X'$ that
contracts all the divisors $E_{1,i}$ using the projection onto
$\PP^{n-2}$. There is an induced rational map $\psi:X\dra X'$ that
is an isomorphism in codimension $1$. Let $E'_1=\psi(E_1)$. Then
$E'_1\cong Y$. In fact, the rational map $X\dra Y$ is resolved by
this flip and induces a regular map $X'\ra Y$ that is a
$\PP^1$-bundle, with $E'_1$ as a section. If $D$ is a divisor on
$X$, let $D'=\psi(D)$. Using geometric arguments, one checks that
on $E'_1\cong Y$ one has $D'_{|E'_1}=\tilde{D}$ when $D=H$,
$H-E_1$, $E_i$, for $i=2,\ldots,r$. Hence, the formula holds in
general by linearity. Then $r$ is the composition of the
isomorphism $\H^0(X,D)\cong\H^0(X',D')$ with the restriction map
$\H^0(X,D')\ra\H^0(E'_1,D'_{|E'_1})$.
\end{rmk}

\begin{notn}
Let $q=\tilde C\cap E_1$.
Obviously, $q\in\pi_1(C)$.
Let $Y'=\Bl_r\PP^{n-1}$ be the blow-up of $Y$ at $q$ and let $E_q$ be the
exceptional divisor.
\end{notn}

%

\begin{lemma}\label{two_cases}
Let $E$ be a minimal divisor on $X$ of $H$-degree $k$. Then
$E\cdot(l-e_1)$ is either $0$ or $1$. In the first case, $\tilde
E$ is a minimal divisor on $Y$. In the second case, the divisor
$E'=\tilde{E}-(k-1)E_q$ is minimal on $Y'$, except when $k=1$. In
the latter case, one has:
\begin{equation}\label{special}
E=H-\sum_{i\in I}E_i, \quad
\tilde{E}=\sum_{i\in I^c}\ov E_{i},
\quad I\subset\{2,\ldots,r\},\quad |I|=n,\quad |I^c|=\al+1.
\end{equation}
\end{lemma}

\bp In the first case:
\begin{equation}\label{originalE_case1}
E=kH-kE_1-k\sum_{i\in I}E_i -(k-1)\sum_{i\in I^c}E_i,
\end{equation}
where $I\subset\{2,\ldots,r\}$,
$|I|=n+1-2k$, and
\begin{equation}\label{E_case1}
\tilde{E}=k\ov H-k\sum_{i\in I}\ov E_i -(k-1)\sum_{i\in I^c}\ov E_i.
\end{equation}

In the second case:
$$
E=kH-(k-1)E_1-k\sum_{i\in I}E_i -(k-1)\sum_{i\in I^c}E_i,$$
where $I\subset\{2,\ldots,r\}$, $|I|=n+2-2k$, and
\begin{equation}\label{E_case2}
\quad\tilde{E}=(k-1)\ov H-(k-1)\sum_{i\in I}\ov E_i -
(k-2)\sum_{i\in I^c}\ov E_i.
\end{equation}

Let $s\in\H^0(Y,\tilde{E})$ be the image of the section $x_E$ via
the map $r$ of (\ref{res2}). Let $Z$ be the zero-locus of $s$. By
Lemma \ref{minimal}, the divisor $E$ has multiplicity $k-1$
along~$\tilde{C}$. Therefore,
$$\mult_q Z=\mult_q E\cap E_1\geq
\mult_q E\geq\mult_{\tilde{C}} E=k-1.
$$
It follows that the image of $r$ is in each case contained
in the push-forward of the linear system $|E'|$ on $Y'$. Except in
Case 2) when $k=1$, $E'$ is minimal on~$Y'$. \ep

We prove that $\H^0(X,D)$ is generated by distinguished sections
for any effective divisor~$D$.

\begin{claim}
We may assume that $0<m_1\le m_2\le\ldots\le m_r$.
\end{claim}

\bp Indeed, if $m_i\leq0$ for some $i$, then
$\H^0(X,D)\cong\H^0(X,D_0)$, where $D_0=D+m_iE_i$ is a divisor on
$\Bl_{r-1}\PP^n$. The ring $\Cox(\Bl_{r-1}\PP^n)$ is minuscule:
this follows by Remark \ref{case r=n+2} if $r=n+3$, and by
induction if $r>n+3$. Hence, $\H^0(X,D_0)$ is generated by
distinguished sections. \ep

\begin{claim}
It suffices to prove that any distinguished section in the image of
$$r:\H^0(X,D)\ra\H^0(Y,\tilde D)$$
can be lifted to a linear combination of distinguished sections.
\end{claim}

\bp Since $\Cox(Y)$ is minuscule by induction and the kernel of
$r$ is $\H^0(X,D-E_1)$, we are then reduced to show that
$\H^0(X,D-E_1)$ is generated by distinguished sections. If $D-E_1$
is effective, we may replace $D$ with $D-E_1$ and repeat the
process. The process stops only when $D-E_1$ is not effective, in
which case $r_{E_1}$ is an isomorphism onto its image. Since for
any effective $D$, one has $D.(l-e_i)=d-m_i\geq0$, for all $i$,
the process must stop. \ep

\begin{notn}
\begin{equation}\label{m}
m= \hbox{max}\left\{\left\lceil\frac{\sum_{i=1}^r
m_i-nd}{\al}\right\rceil,0\right\}
\end{equation}
\end{notn}

\begin{prop}\label{Casem=0}
If $m=0$, then $r$ surjects onto $\H^0(Y,\tilde{D})$ and any
distinguished section $s\in\H^0(Y,\tilde{D})$ can be lifted to  a
distinguished section.
\end{prop}

\bp
The section $s$ is a monomial in the sections corresponding to
minimal divisors on $Y$ and sections $x_{\ov E_i}$,
$i=2,\ldots,r$; hence, it corresponds to a decomposition:
\begin{equation}\label{partition2}
\tilde{D}=S+\sum_{i=2}^rl_i\ov E_i,
\end{equation}
where $l_i\geq0$ and $S$ is a
sum of minimal divisors on $Y$.
Denote
$$\be=d-m_1=D.(l-e_1)\geq0.$$

\begin{lemma}\label{countLemma}
We have $l_i\le\be$ and $\sum\limits_{i=2}^r l_i\geq(\al+1)\be$.
\end{lemma}

\bp For each $k\geq0$, let $a_k\geq0$ be the number of minimal
divisors of $H$-degree $k$ that appear in $S$. Since $\tilde{D}$
and $S$ have the same $H$-degree,
$$m_1=\sum_{k\geq1}ka_k.$$
By counting the number of $\ov E_i$'s in
both sides of (\ref{partition2}), one has
the following formula:

\begin{equation}\label{relation2}
\sum\limits_{i=2}^rl_i=
(\al+1)\be+(nd-\sum\limits_{i=1}^rm_i)+\al(m_1-\sum\limits_{k\geq1}a_k)
\end{equation}
Since $m=0$ and
$$m_1=\sum_{k\geq1}ka_k\geq\sum_{k\geq1}a_k,$$ it
follows that
$$\sum_{i=2}^r l_i\geq(\al+1)\be.$$
Finally,
\begin{equation*}
d-m_i=\tilde{D}.(l-e_i)=S.(l-e_i)+l_i \geq l_i,
\end{equation*}
and therefore
$$l_i\leq(d-m_i)\leq(d-m_1)=\be.$$
\ep

We lift the minimal divisors \eqref{E_case1} on $Y$ to minimal
divisors \eqref{originalE_case1} on $X$ of the same $H$-degree.
Let $D_0$ be the divisor on $X$ equal to the sum of the lifts of
all terms of $S$. Hence, $S=\tilde{D_0}$ and
\begin{equation}\label{deg_diff3}
\hdeg(D)-\hdeg(D_0)=\hdeg(D)-\hdeg(\tilde D)=\beta.
\end{equation}

By Lemma \ref{countLemma} and Lemma \ref{special_lift}, we may
lift $\sum\limits_{i=2}^r l_i\ov E_i$ to an effective divisor
$D_1$ on $X$, with $\hdeg(D_1)=\be$. Let $D'=D_0+D_1$. Then $D'$
has the same $H$-degree as $D$. Since $\tilde{D'}=\tilde{D}$, it
follows from Lemma \ref{trick} that $D=D'$. By construction, there
is a distinguished section $t$ in $\H^0(X,D)$ such that $r(t)=s$.
\ep

\begin{lemma}\label{special_lift}
Consider the divisor $\sum\limits_{i=2}^rl_i\ov E_i$ on $Y$ and
assume that $$l_i\le\be\quad (i=2,\ldots, r), \quad
\sum\limits_{i=2}^r l_i\geq(\al+1)\be.$$

Then we may lift $\sum\limits_{i=2}^rl_i\ov E_i$ to an effective
divisor $D_1$ on $X$ with $\hdeg(D_1)=\be$. Moreover, there is a
distinguished section $t\in\H^0(X,D_1)$ such that
$r(t)=\prod_{i=2}^r x_{\ov E_i}^{l_i}$.
\end{lemma}

\bp For all $i=2,\ldots,r$, we may write $l_i=l'_i+l''_i$, for
some $l'_i,l''_i\geq0$, such that $0\leq l'_i\leq\be$ and
$\sum\limits_{i=2}^rl'_i=(\al+1)\be$. By partitioning
$\sum\limits_{i=2}^rl'_iE_i$ into a sum of $(\al+1)$-tuples of the
form $E_{i_1}+\ldots+E_{i_{\al+1}}$ (the precise procedure for the
partitioning is explained in the proof of Lemma \ref{effective}),
we may lift $\sum\limits_{i=2}^rl'_i\ov E_i$ using (\ref{special})
to a divisor $D'_1$ on $X$ which is a sum of $\be$ ``hyperplane
classes'' $H-\sum E_i$. Hence, $\hdeg(D'_1)=\be$. Moreover, there
is a distinguished section $t'\in\H^0(X,D'_1)$ such that
$r(t')=\prod_{i=2}^r x_{\ov E_i}^{l'_i}$.  Let
$D_1=D'_1+\sum_{i=2}^rl''_i E_i$ and $t=t'\prod_{i=2}^r x_{\ov
E_i}^{l''_i}$. Since $\tilde{E_i}=\ov E_i$ and $r(x_{E_i})=x_{\ov
E_i}$, for all $i=2, \ldots,r$, the Lemma follows. \ep

\begin{prop}
Let $m>0$.  Then the image of $r$ is the push-forward of
$\H^0(Y',\tilde{D}-mE_q)$, and we may lift any distinguished
section $s\in \H^0(Y',\tilde{D}-mE_q)$ to a section $t$ in the
subspace of $\H^0(X,D)$ generated by distinguished sections. By
lift, here we mean that $r(t)=s.x_{E_q}^m$.
\end{prop}

\bp
By Lemma \ref{lower_bound_multip}, the multiplicity of
$D_{|E_1}$ at $q$ is at least $m$. Hence, the map~$r$ has
image in $\H^0(Y',\tilde{D}-mE_q)$.

\begin{lemma}\label{Case12q}
If $E'$ is a minimal divisor on $Y'$ of $H$-degree $k\geq1$, then
the multiplicity at $q$ of a push-forward of $E'$ to $Y$ is
either: 1) $k-1$, or 2) $k$. The push-forward is equal to $\tilde
E$, where $E$ is a minimal divisor on $X$ of $H$-degree $k$ in
Case 1) and $k+1$ in Case~2).
\end{lemma}

\bp
We may lift $E'$
using (\ref{E_case1}) and (\ref{E_case2})
to a minimal divisor $E$ on $X$ by:

\begin{equation*}
1) \quad E=kH-k\sum_{i\in I}E_i -(k-1)\sum_{i\in I^c}E_i,\quad
|I|=n+2-2k, \quad 1\in I
\end{equation*}

\begin{equation*}
2) \quad E=(k+1)H-(k+1)\sum_{i\in I}E_i -k\sum_{i\in I^c}E_i,
\quad |I|=n-2k, \quad 1\in I^c
\end{equation*}

Then $r(x_{E})=x_{E'}x_{E_q}^{k-1}$ in Case 1) and
$r(x_{E})=x_{E'}x_{E_q}^k$ in Case 2).
\ep

Let $S$ be the sum of the minimal divisors $E'$ on $Y'$ whose
sections $x_{E'}$ appear in~$s$. Then:
\begin{equation}\label{partition}
\tilde{D}-mE_q=S+\sum_{i=2}^rl_i\ov E_i+aE_q
\end{equation}
for some integers $l_i,a\geq0$. Hence, the section $s$ in
$\H^0(Y',\tilde{D}-mE_q)$ is of the form~$s'x_{E_q}^{a}$, for $s'$
a section in $\H^0(Y',\tilde{D}-(a+m)E_q)$. So it is enough to
show that we may lift sections $s=s'x_{E_q}^{a}$, with $s'$ a
distinguished section in $\H^0(Y',\tilde{D}-(a+m)E_q)$.

The above lifting $\tilde{E}=E'$ constructs a divisor $D_0$
on $X$ which lifts $S$, i.e.~$S=\tilde D_0$.

\begin{notn}
$$\be=\hdeg(D)-\hdeg(D_0).$$
\end{notn}

If $\be=0$, from Lemma \ref{trick} and $\tilde{E_i}=\ov E_i$ and
$r(x_{E_i})=x_{\ov E_i}$, for all $i=2, \ldots,r$, it follows that
$D=D_0+\sum_{i=2}^rl_iE_i$ and we may lift $s$ to a distinguished
section in $\H^0(X,D)$. For the general case, it is enough to show
that, by eventually rewriting $s$ as a sum of distinguished
sections in $\H^0(Y,\tilde{D})$ corresponding to different
decompositions of $\tilde{D}-mE_q$, we may reduce to the case when
$l_i=l'_i+l''_i$, for some $l'_i,l''_i\geq0$, such that $0\leq
l'_i\leq\be$ and $\sum_{i=2}^rl'_i=(\al+1)\be$. Then we can finish
the proof by using Lemma \ref{special_lift}.

For each $k\geq1$, let $a_k\geq0$, respectively $b_k\geq0$, be the
number of divisors $E'$ as in Case 1), respectively Case 2) of
Lemma \ref{Case12q}, whose sections $x_{E'}$ appear in the
monomial $s$ (taken with multiplicities). One has the following
relations:

\begin{equation}\label{deg_diff2}
0=\hdeg(\tilde{D})-\hdeg(S)=
m_1-\sum_{k\geq1}ka_k-\sum_{k\geq1}kb_k.
\end{equation}

\begin{equation}\label{deg_diff1}
\beta=\hdeg(D)-\hdeg(D_0)=d-\sum_{k\geq1}ka_k-\sum_{k\geq1}(k+1)b_k=
d-m_1-\sum_{k\geq1}b_k.
\end{equation}

Note that by
finding the coefficients of $E_q$ in both sides of the expression
in (\ref{partition}), one has the following relation:
\begin{equation}\label{coef_E_q}
m+a=\sum_{k\geq1}(k-1)a_k+\sum_{k\geq1}kb_k.
\end{equation}

By counting the number of $\ov E_i$'s in both sides of
(\ref{partition}) and using (\ref{deg_diff2}), (\ref{coef_E_q}),
one has:
\begin{equation}\label{relation}
\sum_{i=1}^rm_i-nd=(\al+1)\be+(m+a)\al-\sum_{i=2}^rl_i
\end{equation}

\begin{claim}
We may assume that $a=0$ or $\sum_{k\geq1}b_k=0$.
\end{claim}

\bp Assume $a>0$ and $b_k>0$, for some $k\geq1$. Then the monomial
$s$ contains a section $x_{E'}$, where $E'$ is a minimal divisor
of the form:
$$E'=k\ov H-k\sum_{i\in I}\ov E_i -(k-1)\sum_{i\in I^c}\ov E_i-kE_q,$$
where $I\subset\{2,\ldots,r\}, |I|=n-2k$.  By Lemma \ref{cones},
applied to the divisor $E'+E_q$, we may replace the section
$x_{E'}x_{E_q}$ with a linear combination of sections of the form
$x_{E''}x_{\ov E_j}$, where $j\in\{2,\ldots,r\}$ and
$E''=E'+E_q-\ov E_j$. Then $E''$ is a minimal divisor as in Case
1) of Lemma \ref{Case12q}. Hence, we may replace $s$ with a linear
combination of distinguished sections with smaller $a$ and smaller
$\sum_{k\geq1}b_k$. \ep

\underline{Assume  $\sum b_k=0$.}

Then $\be=d-m_1\geq0$. It follows that $\be\geq l_i\geq0$, for all
$i=2,\ldots,r$. This is because one has from (\ref{partition}):
\begin{equation}\label{d-m_i}
d-m_i=\tilde{D}.(l-e_i)=S.(l-e_i)+l_i\geq l_i
\end{equation}

Hence, $l_i\leq d-m_i\leq\be$, for all $i=2,\ldots,r$.

By definition (\ref{m}), one has $0\leq
m\al-(\sum_{i=1}^rm_i-nd)$. From (\ref{relation}) it follows that
$$(\al+1)\be\leq\sum_{i=2}^rl_i.$$

We are done by Lemma \ref{special_lift}.

\

\underline{Assume $a=0$.}

We show that in this case $\be\geq0$.
By definition (\ref{m}), one has $0\leq
m\al-(\sum_{i=1}^rm_i-nd)<\al$. From (\ref{relation}) it follows
that
$$0\leq\sum_{i=2}^rl_i-(\al+1)\be<\al.$$

It follows that $\be\geq0$. We find $l'_i,l''_i\geq0$ such that
$l_i=l'_i+l''_i$ and $l'_i\leq\be$, for all $i=2,\ldots,r$ and
$\sum_{i=2}^rl'_i=(\al+1)\be$. First, randomly choose $l'_i,l''_i$
with $l_i=l'_i+l''_i$, $l'_i,l''_i\geq0$ and
$\sum_{i=2}^rl'_i=(\al+1)\be$. We show that by eventually
replacing $s$ with a linear combination of distinguished sections
(with smaller $l'_i$), we may reduce to the case when
$l'_i\leq\be$, for all $i$. First take the case when
$i\in\{2,\ldots,r\}$ is such that in $S$ there is no minimal
divisor $E'$ of the form
\begin{equation}\label{E_second_type}
E'=k\ov H-k\ov E_i-k\sum_{j\in I}\ov E_j -(k-1)\sum_{j\in J}\ov
E_j-kE_q
\end{equation}
where $I\subset\{2,\ldots,r\}$, $i\notin I$ and $|I|=n-1-2k$,
$J=\{2,\ldots,r\}\setminus(\{i\}\cup I)$. We claim that
$l_i\leq\be$. Since in each $E'$ appearing in $S$, the divisor
$\ov E_i$ appears with coefficient $-(k-1)$, one has:
\begin{equation*}
d-m_i=\tilde{D}.(l-e_i)=S.(l-e_i)+l_i \geq\sum_{k\geq1}b_k+l_i.
\end{equation*}

It follows that $l_i\leq(d-m_i)-\sum_{k\geq1}b_k\leq\be$.

Assume now that $i\in\{2,\ldots,r\}$ is such that $l'_i>\be$. By
the previous observation, $S$ contains at least one minimal
divisor $E'$ of the form (\ref{E_second_type}). By Lemma
\ref{cones} applied to the divisor $E'+\ov E_i$, we may replace
the section $x_{E'}x_{\ov E_i}$ with a linear combination of
sections of the form $x_{E''}x_{\ov E_j}$, where $j\in J$ and
$E''=E'+E_q-\ov E_j$ is a minimal divisor on $Y'$. Moreover, we
claim that we may choose only indices $j\in J$ with $l'_j<\be$.
Let us call $j\in\{2,\ldots,r\}$ a \emph{good} index if
$l'_j<\be$. We claim that there are at least $k+1$ good indices in
$J$. Clearly, $|J|=r-n+2k-1\geq k+1$. Assume there are at most $k$
good indices in $J$. Then there are at least $(r-n+k-1)=(\al+k+1)$
indices in $J$ that are not good. Since $l'_i>\be$ and $i\notin
J$, it follows that:
$$(\al+1)\be=\sum_{i=2}^rl'_i>(\al+k+1)\be+\be\geq(\al+1)\be$$
which is a contradiction. Hence, the claim follows. By repeating
the process, we end up with $l'_i\leq\be$, for all $i=2,\ldots,r$,
and we are done by Lemma \ref{special_lift}. \ep

\begin{lemma}\label{cones}
Let $X=\Bl_r\PP^n$ be the blow-up of $\PP^n$ in $r\geq n+4$ points
on a rational normal curve $C$ of degree $n$. For any $1\leq k\leq
(n+1)/2$ and any $I\subset\{1,\ldots,r\}$, $|I|=n+1-2k$, let
$$D=kH-k\sum_{i\in I}E_i-(k-1)\sum_{i\in I^c}E_i$$
Then $h^0(X,D)=k+1$. For any $i\in I^c$, the divisor $D-E_i$ is
minimal and, for any choice of $k+1$ indices $i\in I^c$, the
sections $x_{D-E_i}x_{E_i}$ generate $\H^0(X,D)$.
\end{lemma}

\bp
Consider the exact sequence:
\begin{equation}\label{res}
0\ra\H^0(X,D-E_i)\ra\H^0(X,D)\ra\H^0(E_i,D_{|E_i})
\end{equation}

We argue by induction on $n\geq2$.  If $n=2$, then $k=1$ and $D=H-E_j$,
for some $j\in\{1,\ldots,r\}$. Clearly, for any $i\neq j$,
the divisor $H-E_i-E_j$ is minimal.
Since $(H-E_j).E_i=0$, one has $\H^0(E_i,D_{|E_i})\cong\CC$.
For any $l\neq i,j$, the section $x_{H-E_j-E_l}x_{E_l}$ has non-zero
restriction to $E_i$. Hence, the map $\H^0(X,D)\ra\H^0(E_i,D_{|E_i})$
is surjective and $\H^0(X,D)$ is generated by the sections
$x_{H-E_i-E_j}x_{E_i}$ and $x_{H-E_j-E_l}x_{E_l}$.

Assume $n\geq3$. Fix some $i\in I^c$. From Lemma \ref{minimal} the
divisor $E=D-E_i$ is a minimal divisor. Let $Y=\Bl_{r-1}\PP^{n-1}$
be the blow-up of $\PP^{n-1}$ in $r-1$ points corresponding to the
projection from $p_i$ and let $Y'=\Bl_r\PP^{n-1}$ be the blow-up
of $Y$ at the extra point $q$. Then the restriction map in
(\ref{res}) factors through the map
$r_{E_i}:\H^0(X,D)\ra\H^0(Y,\tilde{D})$, where
$$\tilde{D}=(k-1)H-(k-1)\sum_{j\in I}E_j
-(k-2)\sum_{j\in I^c\setminus\{i\}}E_j.$$

Note that by Lemma \ref{lower_bound_multip}, the multiplicity at
$q$ of any divisor in the linear system $|D|$ is at least:
$$\frac{k(n+1-2k)+(k-1)(r+2k-n-1)-nk}{r-n-2}
=k-1-\frac{1}{r-n-2}.$$

Since $r\geq n+4$, the map $r_{E_i}$ has image in $\H^0(Y',D')$,
where $D'=\tilde{D}-(k-1)E_q$. By induction, $\H^0(Y',D')$ has
dimension $k$ and it is generated by any distinct $k$ sections of
the form $x_{E'}x_{E_j}$, where $E'=D'-E_j$ and $j\in
I^c\setminus\{i\}$. On $X$, the divisor $E=D-E_j$ is minimal. By
Lemma \ref{Case12q}, $r_{E_i}(x_E)=x_{E'}x_{E_q}^{k-1}$. Since
$r_{E_i}(x_{E_j})=x_{E_j}$, it follows that
$r_{E_i}(x_Ex_{E_j})=x_{E'}x_{E_j}x_{E_q}^{k-1}$. Hence, the map
$r_{E_i}$ has image $\H^0(Y',D')x_{E_q}^{k-1}$. Therefore,
$\H^0(X,D)$ has dimension $k+1$ and it is generated by any $k+1$
sections of the form  $x_Ex_{E_i}$, where $i\in I^c$, $E=D-E_i$.
\ep

\begin{thm}\label{n=2}  Let $r\geq5$ and let $X=\Bl_r\PP^2$ be the blow-up
of $\PP^2$ at $r$ distinct points $p_1,\ldots,p_r$ that lie on an
irreducible conic. Then $\Cox(X)$ is minuscule.
\end{thm}

\bp Let $C$ be the proper transform on $X$ of the conic in $\PP^2$
that contains the points $p_1,\ldots,p_r$. Then
$C=2H-\sum_{i=1}^rE_i$. For any $i,j\in\{1,\ldots,r\}$ with $i\neq
j$, let $L_{i,j}$ be the proper transform on $X$ of the line that
passes through $p_i$, $p_j$. The classes $C$ and $L_{i,j}$ are the
minimal divisors on $X$. Let $x_C$, resp. $x_{L_{i,j}}$, be the
corresponding sections. A distinguished section on $X$ is a
monomial in $x_C$, $x_{L_{i,j}}$ and $x_{E_i}$, for all $i,j$.

We prove by induction on $r$ that $\Cox(\Bl_r\PP^2)$ is generated
by distinguished sections. The case $r=5$ was proved in \cite{BP}.
Assume $r\geq6$.

Let $D$ be an effective divisor \eqref{divisorD} on $X$. If
$m_i=D.E_i\leq0$ for some $i\in\{1,\ldots,r\}$, then
$\H^0(X,D)\cong\H^0(X,D_0)$, where $D_0=D+m_iE_i$ is a divisor on
$\Bl_{r-1}\PP^n$ and $\H^0(X,D_0)$ is generated by distinguished
sections by induction. It follows that $\H^0(X,D)$ is generated by
distinguished sections (obtained by multiplying sections of
$\H^0(X,D_0)$ by $x_{E_i}^{-m_i}$). Hence, we may assume that
$d,m_i>0$ and argue by induction on $d$.

From the exact sequence
$$0\ra\H^0(X,D-C)\ra\H^0(X,D)\ra\H^0(C,D_{|C})$$
it follows that if $D.C=-a<0$, then $\H^0(C,D_{|C})=0$ and
$\H^0(X,D)\cong\H^0(X,D-C)$ is generated by global
sections by induction.

Assume now $D.C=2d-\sum_{i=1}^rm_i\geq0$ and $m_i>0$ for all
$i=1,\ldots r$. Without loss of generality, we may assume $m_1\leq
m_i$, for all $i$. Consider the exact sequence:
$$0\ra\H^0(X,D-E_1)\ra\H^0(X,D)\ra\H^0(E_1,D_{|E_1}).$$

Note $\H^0(E_1,D_{|E_1})=\H^0(\PP^1,\O(m_1))$. For $i=2,\ldots,r$,
let $q_i=L_{1,i}\cap E_1$. Let $x_i\in\H^0(\PP^1,\O(1))$ be the
section vanishing at $q_i$. The divisor $D_{|E_1}$ has
multiplicity at least $m_1+m_i-d$ at $q_i$. Let
$I\subset\{2,\ldots,r\}$ be the set of indices $i$ for which
$m_1+m_i-d\geq0$. It follows that the image of the restriction map
\begin{equation}\label{res_map2}
r:\H^0(X,D)\ra\H^0(E_1,D_{|E_1})
\end{equation}
lies in the subspace
$$V=\prod_{i\in I}x_i^{m_1+m_i-d}\H^0(\PP^1,\O(e))
\subset\H^0(\PP^1,\O(m_1)),$$
where
\begin{equation}\label{e}
e=m_1-\sum_{i\in I}(m_1+m_i-d)
\end{equation}

We claim that one may lift any section in $V$ to a section in
$\H^0(X,D)$ that is generated by distinguished sections.
Then we are reduced to show that $\H^0(X,D-E_1)$ is generated by
distinguished sections. If $D-E_1$ is not effective, we are done; if not,
we replace $D$ with $D-E_1$ and repeat the process until either $D-E_1$ is
not effective or $D.E_i\leq0$.

Clearly, $\H^0(\PP^1,\O(e))$ is generated by sections
$\prod_{i=2}^r x_i^{k_i}$, where $k_i\ge0$ and $\sum k_i=e$ (of
course, we may assume that, for example, $k_4=k_5=\ldots=0$). Note
that $r(x_{L_{1,j}})=x_j$, for all $j=2,\ldots,r$. Consider the
following divisor on $X$:
$$D_0=\sum_{i\in I^c}k_iL_{1,i}+\sum_{i\in I}(k_i+m_1+m_i-d) L_{1,i}$$
$$=m_1H-m_1E_1-\sum_{i\in
I^c}k_iE_i-\sum_{i\in I}(k_i+m_1+m_i-d)E_i.$$

The restriction map $r$ maps the section
$$t'=\prod_{i=2}^r
x_{L_{1,i}}^{k_i}\prod_{i\in I}x_{L_{1,i}}^{m_1+m_i-d}\in
\H^0(X,D_0)$$ to the section
$$s=\prod_{i=2}^r x_i^{k_i}\prod_{i\in
I}x_i^{m_1+m_i-d}\in \H^0(E_1,{D_0}_{|E_1})=\H^0(\PP^1,\O(m_1)).$$
Consider
$$D-D_0=(d-m_1)H-\sum_{i\in I^c} (m_i-k_i)E_i-\sum_{i\in
I}(d-m_1-k_i)E_i.$$
Since
$$d\ge m_1,\quad
m_i\ge m_1\ge e\ge k_i,\quad d-m_1\ge m_1\ge e\ge k_i$$ and using
(\ref{e}) one has:
$$\sum_{i\in I^c}(m_i-k_i)+\sum_{i\in I}(d-m_1-k_i)=
\sum_{i\in I^c}m_i+(d-m_1)|I|-e=\sum_{i=1}^r
m_i-2m_1\leq2(d-m_1).$$ It follows from Lemma \ref{effective} that
$D-D_0$ is an effective divisor on $X$. Since $(D-D_0).E_1=0$, the
space $\H^0(X,D-D_0)$ is generated by distinguished sections by
induction. Let $t''\in\H^0(X,D-D_0)$ be any distinguished section
not zero on $E_1$. Then $t't''$ is a distinguished section in
$\H^0(X,D)$ that maps to $s$. \ep

\begin{lemma}\label{effective}
Let $X$ be the blow-up of $\PP^n$ in any $r$ distinct points. Let
$D=dH-\sum_{i=1}^rm_iE_i$, $d,m_i\ge0$, be a divisor class with
$\sum_{i=1}^rm_i\leq nd$ and $d\geq m_i$, for all $i=1,\ldots r$.
Then $D$ is an effective divisor.
\end{lemma}

\bp We claim that $D$ is an effective combination of (effective)
classes $H-(E_{i_1}+\ldots+E_{i_l})$, for
$i_1,\ldots,i_l\in\{1,\ldots,r\}$ and $0\leq l\leq n$.
Consider the table with $n$ rows and $d$ columns filled with $E_i$'s
in the following way. Start in the upper left corner and write $m_1$
$E_1$'s in the first row. Then write $m_2$ $E_2$'s passing to
the second row if necessary, and so on.
Fill the remaining entries with zeros.
In the following example $n=3$ and $D=5H-3E_1-3E_2-2E_3-5E_4-E_5$:
$$\begin{matrix}
E_1&E_1&E_1&E_2&E_2\cr
E_2&E_3&E_3&E_4&E_4\cr
E_4&E_4&E_4&E_5&0\cr
\end{matrix}
$$
Our conditions guarantee that all entries of a given column are different.
Therefore $D$ is the sum of classes $H-(E_{i_1}+\ldots+E_{i_l})$,
one for each column,
where $E_{i_1},\ldots,E_{i_l}$ are entries of the column.
In the example above,
$$D=(H-E_1-E_2-E_4)+(H-E_1-E_3-E_4)+\ldots+(H-E_2-E_4).$$
\ep

\section{Proof of Theorem~\ref{Thm_intro}}\label{generators}
By \cite{Mukai01}, there is
an isomorphism $\phi:S^G\ra\Cox(X)$ where $X$ is the blow-up of
$\PP^n$ in $n+3$ points $p_1,\ldots,p_{n+3}$ in general position.
By Theorem~\ref{MinimalDivisors}, the ring $\Cox(X)$ is generated
by the sections $x_{E_i}$, for each exceptional divisor $E_i$,
$i=1,\ldots,n+3$, and the sections $x_E$, corresponding to the
minimal divisors
\begin{equation}\label{exc}
E=kH-k\sum_{i\in I}E_i -(k-1)\sum_{i\in I^c}E_i
\end{equation}
for each subset $I\subset\{1,\ldots,n+3\}$, $|I|=n+2-2k$, $1\leq
k\leq1+n/2$. Then $|I^c|=2k+1$. Note that if $k=0$ in (\ref{exc}),
then $E=E_i$.

The polynomials $F_I$ in (\ref{invariants}) are clearly
invariant (just use the rule of differentiating a determinant).
We claim that, for all $0\leq k\leq1+n/2$ one has
$\phi(F_{I^c})=x_E$, where $E$ is as in (\ref{exc}). It is clear from
\cite{Mukai01} that $\phi(x_i)=x_{E_i}$.  Following \cite{Mukai01},
if $F_0=\ldots=F_n=0$ are $n+1$ linear equations (in $t_1,\ldots,t_r$)
that cut $G$ in $\bG^{n+3}$, let $J_0,\ldots J_n$ be the polynomials in $S$
given by $J_i=F_i(y_1/x_1,\ldots,y_r/x_r)x_1\ldots x_r$.
Then sections of the divisor
$D=dH-\sum_{i=1}^{n+3}m_iE_i$ on $X$, for $d,m_i\geq0$,
correspond by $\phi$ to an invariant polynomial of the form
$$Q=\frac{P(J_0,\ldots,J_n)}{\prod_{i=1}^{n+3}x_i^{m_i}}$$ where
$P(z_0,\ldots,z_n)$ is a homogeneous polynomial of degree $d$ in
variables $z_0,\ldots,z_n$, such that $P(J_0,\ldots,J_n)$ is
divisible by $\prod_{i=1}^{n+3}x_i^{m_i}$.  If we let $\deg_x(Q)$, resp.
$\deg_y(Q)$, to be the degree of $Q$ in the $x_i$'s, resp. in the $y_i$'s,
then
\begin{equation}\label{Q}
\deg_y(Q)=d,\quad\deg_x(Q)=(n+2)d-\sum_{i=1}^{n+3}m_i,
\quad\deg(D)=\deg_x(Q)-\deg_y(Q).
\end{equation}

Hence, $\phi(F_{I^c})$ is a section in $\H^0(X,D)$, where $D$ is a divisor
with $d=k$ and $\deg(D)=1$. To show that $D=E$, consider the following
action of the torus $\bG_m^r$ on $S$: $(\la_1,\ldots,\la_r)\in\bG_m^r$ acts
by $x_i\mapsto\la_i x_i$, $y_i\mapsto\la_i y_i$. The action of  $\bG_m^r$
on $S$ is compatible with the action of  $\bG_a^r$ on $S$. Hence, there is
an induced action of $\bG_m^r$ on $S^G$. Since
$(\la_1,\ldots,\la_r)\in\bG_m^r$ maps $J_i$ onto
$\la_1\ldots\la_r J_i$, it follows that $Q$ is mapped to
$\prod_{i=1}^r\la_i^{d-m_i}$. Since $F_{I^c}$ is mapped to
$\prod_{i\in I^c}\la_i$, it follows that $D=E$.
\qed

\end{document}